\documentclass[review]{elsarticle}
\usepackage{amssymb,latexsym,amsmath}
\DeclareMathOperator{\sech}{sech}
\usepackage{graphicx,times}
\usepackage[margin=1.0in]{geometry}
\newtheorem{Ex}{Example}[section]
\numberwithin{equation}{section}
\journal{Journal of \LaTeX\ Templates}

\begin{document}

\begin{frontmatter}

\title{ Modified DJ method: Application to Boussinesq equation}

\author{Jayvant Patade$^*$\fnref{myfootnote} }
\ead{dr.jayvantpatade@gmail.com}
\author{Sachin Bhalekar$^2$ }
\ead{ sachin.math@yahoo.co.in}
\address{$^{*1}$ Department of Mathematics, Jaysingpur College, Jaysingpur, Kolhapur, India - 416101\\$^2$ School of Mathematics and Statistics, University of Hyderabad, India - 500046}


\cortext[mycorrespondingauthor]{Corresponding author}
\begin{abstract}
 In this paper we present a  modification of  DJ Method [J. Math. Anal. Appl. 316 (2006), 753-763] to solve the nonlinear equations more efficiently. It is observed that the modified DJ method is faster and hence it has accelerated convergence rate as compared to the old one. We use this new method to find the analytical solutions of Boussinesq equation. The reported results  are compared with the exact solutions. Further, we compare the absolute error in our solution with those in other iterative methods. It is observed that the presented method is simple and generates more accurate solutions as compared with other methods.
\end{abstract}
\end{frontmatter}

\section{Introduction}
\par Boussinesq equation introduced by French mathematician Joseph Boussinesq has the form
\begin{equation}
u_{tt} + p u_{xx} + q \left( u^2\right)_{xx} + r u_{xxxx} = 0, \label{1}
\end{equation}
where p, q and r are constants. The Boussinesq equation have several applications in the real world. This equation play an important role in modeling various phenomena such as long waves in shallow water \cite{BE}, one dimensional nonlinear lattices waves \cite{lattices}, vibration in a nonlinear string \cite{string}, electromagnetic waves in dielectric materials \cite{dielectric} and so on. Many researchers have been used analytical methods to solve Boussinesq equation such as variational iteration method \cite{VIM}, modified variational iteration method \cite{MVIM1,MVIM2},  Adomian decomposition method and homotopy perturbation method \cite{Mohyud,ADM1}. Recently Malek et al. \cite{Malek} have used potential symmetries method to solve Boussinesq equation.

\par The Daftardar-Gejji and Jafari Method (DJM) \cite{GEJJI1} is a simpler and more efficient technique  used to solve various equations such as  fractional differential equations \cite{GEJJI3}, partial differential equations \cite{GEJJI4}, boundary value problems \cite{GEJJI5}, evolution equations \cite{GEJJI6}, system of nonlinear functional equations \cite{GEJJI8}, algebraic equations \cite{Noor} and so on. The method is successfully employed to solve Newell-Whitehead-Segel equation \cite{pb1}, Fisher’s equation \cite{pb8}, Ambartsumian equation \cite{pb5}, fractional-order logistic equation \cite{GEJJI7} and some nonlinear  dynamical systems \cite{GEJJI9,pb4}. In \cite{pb2,pb3} we provided the series solutions of pantograph equation in terms of new special functions. Recently DJM has been used to generate a new numerical methods \cite{pb6,pb7} for solving  differential equations.

\par In this manuscript we consider the solution of Boussinesq equation by using modified DJM. We organize the paper as follows:
\par The DJM is described briefly in Section \ref{djm} and the modified DJM is described in section \ref{mdjm}. A general technique used to solve Boussinesq equation  using modified DJM is described in Section \ref{appl}. Section \ref{Ex} deals with  illustrative examples and the conclusions are summarized in Section \ref{concl}.

\section{Daftardar-Gejji and Jafari Method}\label{djm}
In this section we describe the Daftardar-Gejji and Jafari Method (DJM) \cite{GEJJI1}, which is useful for solving the nonlinear equations of the form
\begin{equation}
u= f + L(u) + N(u),\label{2.1}
\end{equation}
where $f$ is a source term, $L$ and $N$ are linear and nonlinear operators respectively. It is assumed that the DJM  solution for the  Eq.(\ref{2.1}) has the form:
\begin{equation}
u= \sum_{i=0}^\infty u_i. \label{2.2}
\end{equation}
The convergence of series (\ref{2.2}) is proved in 	\cite{CONV}.
\par Since $L$ is linear
\begin{equation}
L\left(\sum_{i=0}^\infty u_i\right) = \sum_{i=0}^\infty L(u_i).\label{2.3}
\end{equation}
The nonlinear operator $N$ in Eq.(\ref{2.1}) is decomposed by Daftardar-Gejji and Jafari as bellow:
\begin{eqnarray}
N\left(\sum_{i=0}^\infty u_i\right) &=& N(u_0) + \sum_{i=1}^\infty \left\{N\left(\sum_{j=0}^i u_j\right)- N\left(\sum_{j=0}^{i-1} u_j\right)\right\} \nonumber\\
&=& \sum_{i=0}^\infty G_i,\label{2.4}
\end{eqnarray}
\\
where $G_0=N(u_0)$ and $G_i = \left\{N\left(\sum_{j=0}^i u_j\right)- N\left(\sum_{j=0}^{i-1} u_j\right)\right\}$, $i\geq 1$.\\
\\
Using equations (\ref{2.2}), (\ref{2.3}) and (\ref{2.4}) in Eq.(\ref{2.1}),  we get
\begin{equation}
\sum_{i=0}^\infty u_i= f + \sum_{i=0}^\infty L(u_i) + \sum_{i=0}^\infty G_i.\label{2.5}
\end{equation}
From Eq.(\ref{2.5}), the DJM series terms are generated as bellow:
\begin{eqnarray}
u_0 &=& f,\nonumber\\
u_{m+1} &=& L(u_m) + G_m,\quad m=0,1,2, \cdots.\label{2.6}
\end{eqnarray}
In practice, we take the approximation
\begin{eqnarray}
u = \sum_{i=0}^{k-1} u_i
\end{eqnarray}
for suitable integer $k$.\\
The convergence results for DJM are described in \cite{CONV}.
\section{Modified Daftardar-Gejji and Jafari Method}\label{mdjm}
\par In  \cite{Wazwaz} Wazwaz proposed a  modification in ADM to generate a rapidly converging solution series. Using same technique we modify the Daftardar-Gejji and Jafari method as follows:\\
We assume that
\begin{equation}
f = f_1 + f_2.\label{3.1}
\end{equation}
Then Eq.(\ref{2.1}) can be written as
\begin{equation}
u=  f_1 + f_2 + L(u) + N(u).\label{3.2}
\end{equation}
The modified DJM is described  as:
\begin{eqnarray}
u_0 &=& f_1,\nonumber\\
u_1 &=& f_2 + L(u_0) + G_0,\nonumber\\
u_{m+1} &=& L(u_m) + G_m,\quad m=1,2,3, \cdots.\label{3.3}
\end{eqnarray}
It is obvious that the simpler form of $u_0$ in MDJM result in the reduction of computations and accelerates the convergence rate.
The convergence results for MDJM are similar as that of DJM \cite{CONV}

\section{Applications}\label{appl}

A bidirectional solitary wave solution of Eq.(\ref{1}) discussed in  \cite{Clarkson} and is given by

\begin{equation}
u(x,1) =- \frac{3(\alpha^2+p)^{\frac{1}{2}}}{2q} \sech^2\left[\frac{1}{2}\left(\frac{\alpha^2+p}{-r}\right)^{\frac{1}{2}}(x\pm\alpha t) + \beta\right],
\end{equation}
where $\alpha$ and $\beta$ are constants and the initial conditions are
\begin{eqnarray}
u(x,0) &=& \frac{3(\alpha^2+p)^{\frac{1}{2}}}{2q} \sech^2\left[\frac{1}{2}\left(\frac{\alpha^2+p}{-r}\right)^{\frac{1}{2}}x +\beta\right],\\
u_t(x,0) &=& \mp\frac{3\alpha(\alpha^2+p)^{\frac{3}{2}}}{2q(-r)^{\frac{1}{2}}} \sech^2\left[\frac{1}{2}\left(\frac{\alpha^2+p}{-r}\right)^{\frac{1}{2}}x +\beta\right]\nonumber\\
&&\tanh\left[\frac{1}{2}\left(\frac{\alpha^2+p}{-r}\right)^{\frac{1}{2}}x +\beta\right].
\end{eqnarray}
The equivalent integral equation of (\ref{1}) is
\begin{eqnarray}
u(x,t) &=& \frac{3(\alpha^2+p)^{\frac{1}{2}}}{2q} \sech^2\left[\frac{1}{2}\left(\frac{\alpha^2+p}{-r}\right)^{\frac{1}{2}}x +\beta\right]\nonumber\\
&&\mp\frac{3\alpha(\alpha^2+p)^{\frac{3}{2}}}{2q(-r)^{\frac{1}{2}}} \sech^2\left[\frac{1}{2}\left(\frac{\alpha^2+p}{-r}\right)^{\frac{1}{2}}x +\beta\right]\nonumber\\
&&\tanh\left[\frac{1}{2}\left(\frac{\alpha^2+p}{-r}\right)^{\frac{1}{2}}x +\beta\right] t\nonumber\\
&& - \int_0^t(u_{xx} + u_{xxxx}) dx - \int_0^t  \left(u^2\right)_{xx} dx. \label{4.1}
\end{eqnarray}
This Eq.(\ref{4.1}) is of the form Eq.(\ref{3.2}). Using Eq.(\ref{3.3}), the modified DJM series terms are generated as bellow:
\begin{equation}
u_0(x,t) = - \frac{3(\alpha^2+p)^{\frac{1}{2}}}{2q} \sech^2\left[\frac{1}{2}\left(\frac{\alpha^2+p}{-r}\right)^{\frac{1}{2}}x +\beta\right],
\end{equation}
\begin{eqnarray}
u_1(x,t) &=& \mp\frac{3\alpha(\alpha^2+p)^{\frac{3}{2}}}{2q(-r)^{\frac{1}{2}}} \sech^2\left[\frac{1}{2}\left(\frac{\alpha^2+p}{-r}\right)^{\frac{1}{2}}(x\pm\alpha t +\beta\right]\nonumber\\
&&\tanh\left[\frac{1}{2}\left(\frac{\alpha^2+p}{-r}\right)^{\frac{1}{2}}(x\pm\alpha t +\beta\right]t\nonumber\\
&& - \int_0^t\left((u_0(x,t))_{xx} + (u_0(x,t))_{xxxx}\right) dx\\
&& - \int_0^t  (u_0^2(x,t))_{xx} dx,\nonumber\\
u_{n+1}(x,t)&=& -\int_0^t\left(\left(\sum_{i=0}^{n} u_i(x,t)\right)_{xx} + \left(\sum_{i=0}^{n}u_i(x,t)\right)_{xxxx} \right)dx\nonumber\\
&& - \int_0^t \left(\sum_{i=0}^{n}u_i(x,t)\right)^2_{xx} dx \nonumber\\
&& + \int_0^t \left(\sum_{i=0}^{n-1}u_i(x,t)\right)^2_{xx} dx,\quad n=1,2,3,\cdots.
\end{eqnarray}

\section{Illustrative examples}\label{Ex}
Besides equation Eq.(\ref{1}), there are few more PDEs which are called  Boussinesq equation. In this section, we solve such equations using MDJM.
\begin{Ex}\label{ex1}
	Consider the  Boussinesq equation \cite{Mohyud}
	\begin{equation}
	u_{tt} - u_{xx} + 3 \left( u^2\right)_{xx} + u_{xxxx} = 0, \label{5.1}
	\end{equation}
	with initial condition
	\begin{eqnarray}
	u(x,0) &=& \frac{c}{2} \sech^2\left[\frac{\sqrt{c}}{2}(x+1)\right],\\
	u_t(x,0) &=& -\frac{c^\frac{5}{2}}{2} \sech^2\left[\frac{\sqrt{c}}{2}(x+1)\right]\tanh\left[\frac{\sqrt{c}}{2}(x+1)\right] .
	\end{eqnarray}
\end{Ex}
The equivalent integral equation is
\begin{eqnarray}
u(x,t) &=& \frac{c}{2} \sech^2\left[\frac{\sqrt{c}}{2}(x+1)\right] - \frac{c^\frac{5}{2}}{2} \sech^2\left[\frac{\sqrt{c}}{2}(x+1)\right]\tanh\left[\frac{\sqrt{c}}{2}(x+1)\right] t \nonumber\\
&& + \int_0^t(u_{xx} - u_{xxxx}) dx - 3 \int_0^t  u^2_{xx} dx. \label{5.2}
\end{eqnarray}
This Eq.(\ref{5.2}) is of the form Eq.(\ref{3.2}). Using Eq.(\ref{3.3}), the MDJM series terms are generated as bellow:
\begin{eqnarray}
u_0(x,t) &=& \frac{c}{2} \sech^2[\frac{\sqrt{c}}{2}(x+1)],\\
u_1(x,t) &=& -\frac{1}{8} c^2 t^2 \sech^4\left[\frac{1}{2} \sqrt{c} (1+x)\right]-\frac{1}{4} c^3 t^2 \sech^6\left[\frac{1}{2}\sqrt{c} (1+x)\right]\nonumber\\
&&-\frac{1}{4} c^{5/2} t \sech^2\left[\frac{1}{2} \sqrt{c} (1+x)\right] \tanh\left[\frac{1}{2} \sqrt{c} (1+x)\right]+\cdots,\\
u_2(x,t) &=& \frac{1}{48} c^3 t^4 \sech^6\left[\frac{1}{2} \sqrt{c} (1+x)\right]+\frac{13}{192} c^4 t^4 \sech^8\left[\frac{1}{2}\sqrt{c} (1+x)\right]\nonumber\\
&& -\frac{17}{192} c^5 t^4 \sech^{10}\left[\frac{1}{2} \sqrt{c} (1+x)\right]-\cdots.,\\
u_3(x,t) &=& \frac{7}{64} c^4 t^4 \sech^8\left[\frac{1}{2} \sqrt{c} (1+x)\right]\nonumber\\
&&-\frac{17 c^4 t^6 \sech^8\left[\frac{1}{2}\sqrt{c} (1+x)\right]}{5760}\nonumber\\
&&+\frac{47}{64} c^5 t^4 \sech^{10}\left[\frac{1}{2} \sqrt{c} (1+x)\right]\nonumber\\
&&-\frac{77 c^5 t^6 \sech^{10}\left[\frac{1}{2}\sqrt{c} (1+x)\right]}{1920}+\cdots.
\end{eqnarray}
and so on.

The exact solution of Eq.(\ref{5.1}) is
\begin{equation}
u(x,t) = \frac{c}{2} \sech^2\left[\frac{\sqrt{c}}{2}x+\frac{\sqrt{c}}{2}\sqrt{1+ct}\right]
\end{equation}
\\
We compare 4-term solutions for $c=1$ and $c=2$ in Fig.1 and Fig.2, where MDJM solution and exact solution are shown by red and  green  colors respectively. From these figures, it can be observed that the  modified DJM solution is well in agreement to exact solution.

\begin{tabular}{c}
	\includegraphics[scale=1]{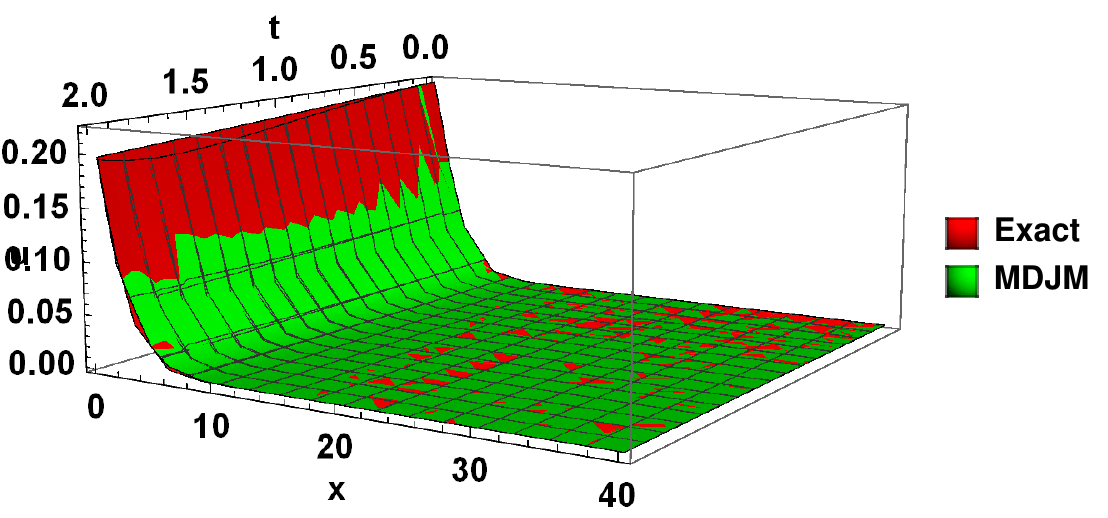} \\
	Fig.1: Comparison of solutions of Eq.(\ref{5.1}) for $c=1$
\end{tabular}

\begin{tabular}{c}
	\includegraphics[scale=1]{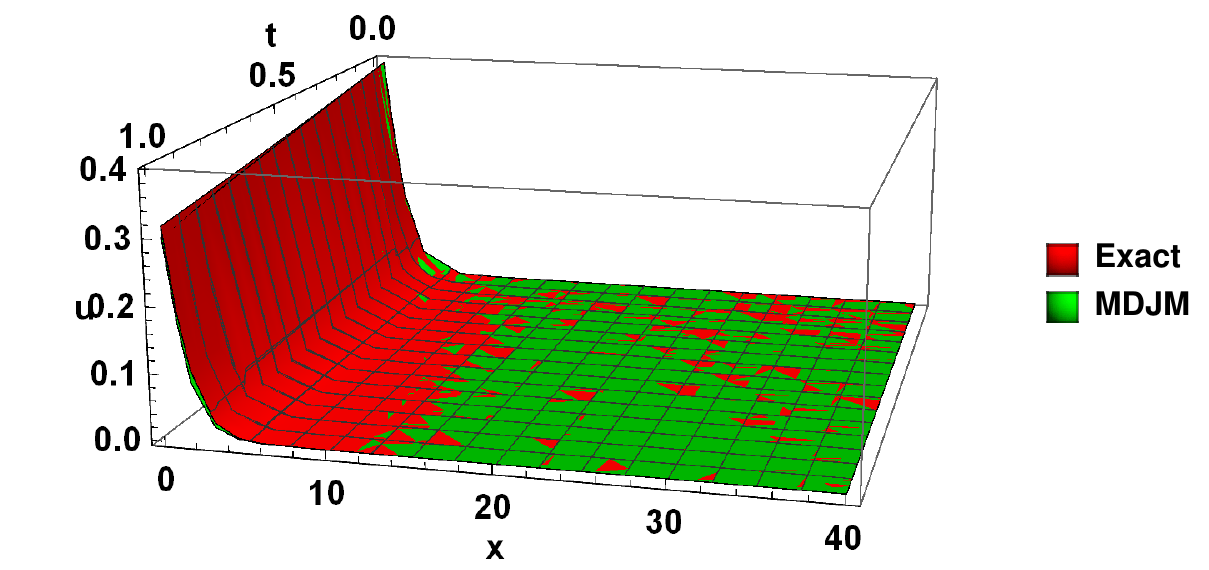} \\
	Fig.2: Comparison of solutions of Eq.(\ref{5.1}) for $c=2$
\end{tabular}

The 6-term ADM and HPM solutions of (\ref{5.1}) are given in \cite{Mohyud} as\\
\begin{eqnarray}
u(x,t) &=& \frac{1}{2} \sech^2\left[\frac{1}{2} \sqrt{c} (1+x)\right]\nonumber\\
&&-\frac{1}{8}(-1+c) c^2 t^2\left(-2+\cosh[\sqrt{c}(x+1)]\sech^4\left[\frac{1}{2} \sqrt{c} (1+x)\right]\right.\nonumber\\
&&-\frac{1}{1280}\left((-1+c)^2 c^5 t^6 (-140 + 157 \cosh[\sqrt{c}(x+1)])\right) -26 \cosh[2\sqrt{c}(x+1)] \nonumber\\
&&+ \cosh[3\sqrt{c}(x+1)]\sech^{10}\left[\frac{\sqrt{c}}{2}(1+x)\right]+\cdots\label{5.4}
\end{eqnarray}
\\
and
\begin{eqnarray}
u(x,t) &=& \frac{1}{2} \sech^2\left[\frac{1}{2} \sqrt{c} (1+x)\right]\nonumber\\
&&-\frac{1}{8}(-1+c) c^2 t^2\left(-2+\cosh[\sqrt{c}(x+1)]\sech^4\left[\frac{1}{2} \sqrt{c} (1+x)\right]\right.\nonumber\\
&&-\frac{1}{30720}\left((-1+c) c^3 t^4 \left(-475 - 6125c - 9480 c^2 - 3360c^2t^2 +\cdots\right.\right.\label{5.5}
\end{eqnarray}
respectively. We compare the errors in all these solutions in table 1 and 2.
\\
\\
Table 1: Absolute error in the 4-term MDJM solution and exact solution for $c=1$.\\
\\
\begin{tabular}{|c|c|c|c|c|c|}
	\hline
	$t/x$ & $20$ & $25$ & $30$ & $35$ & $40$\\
	\hline
	$0.1$ & $3.58376\times10^{-12}$  & $2.41472\times10^{-14}$ & $1.62702\times10^{-16}$ & $1.09628\times10^{-18}$ & $7.38668\times10^{-21}$ \\
	\hline
	$0.2$ & $1.36005\times10^{-11}$  & $9.16392\times10^{-14}$ & $6.1746\times10^{-16}$ & $4.16041\times10^{-18}$ & $2.80326\times10^{-20}$ \\
	\hline
	$0.3$ & $2.91257\times10^{-11}$  & $1.96248\times10^{-13}$ & $1.32231\times10^{-16}$ & $8.90963\times10^{-18}$ & $6.00326\times10^{-20}$ \\
	\hline
	$0.4$ & $4.94206\times10^{-11}$  & $3.32993\times10^{-13}$ & $2.24369\times10^{-15}$ & $1.51179\times10^{-17}$ & $1.01863\times10^{-19}$ \\
	\hline
	$0.5$ & $7.38844\times10^{-11}$  & $4.97829\times10^{-13}$ & $3.35435\times10^{-15}$ & $2.26014\times10^{-17}$ & $1.52287\times10^{-19}$ \\
	\hline
	$0.6$ & $1.02022\times10^{-10}$  & $6.8742\times10^{-13}$ & $4.6318\times10^{-15}$ & $3.12088\times10^{-17}$ & $2.10283\times10^{-19}$ \\
	\hline
	$0.7$ & $1.33421\times10^{-10}$  & $8.98981\times10^{-13}$ & $6.05728\times10^{-15}$ & $4.08137\times10^{-17}$ & $2.7500\times10^{-19}$ \\
	\hline
	$0.8$ & $1.67731\times10^{-10}$  & $1.13017\times10^{-12}$ & $7.61499\times10^{-15}$ & $5.13094\times10^{-17}$ & $3.4572\times10^{-19}$ \\
	\hline
	$0.9$ & $2.04658\times10^{-10}$  & $1.37898\times10^{-12}$ & $9.29146\times10^{-15}$ & $6.26054\times10^{-17}$ & $4.21832\times10^{-19}$ \\
	\hline
	$1$ & $2.43946\times10^{-10}$  & $1.64369\times10^{-12}$ & $1.10751\times10^{-14}$ & $7.46236\times10^{-17}$ & $5.0281\times10^{-19}$ \\
	\hline
\end{tabular}
\\
\\
Table 2: Absolute error in the 4-term MDJM solution and exact solution for $c=2$.\\
\\
\begin{tabular}{|c|c|c|c|c|c|}
	\hline
	$t/x$ & $20$ & $25$ & $30$ & $35$ & $40$\\
	\hline
	$0.1$ & $1.24823\times10^{-14}$  & $1.06015\times10^{-17}$ & $9.00415\times10^{-21}$ & $7.64746\times10^{-24}$ & $6.49518\times10^{-27}$ \\
	\hline
	$0.2$ & $4.58263\times10^{-14}$  & $3.89215\times10^{-17}$ & $3.3057\times10^{-20}$ & $2.80762\times10^{-23}$ & $2.38458\times10^{-26}$ \\
	\hline
	$0.3$ & $9.50542\times10^{-14}$  & $8.0732\times10^{-17}$ & $6.85678\times10^{-20}$ & $5.82364\times10^{-23}$ & $4.94616\times10^{-26}$ \\
	\hline
	$0.4$ & $1.56166\times10^{-13}$  & $1.32635\times10^{-16}$ & $1.12651\times10^{-19}$ & $9.56772\times10^{-23}$ & $8.12611\times10^{-26}$ \\
	\hline
	$0.5$ & $2.25727\times10^{-13}$  & $1.91716\times10^{-16}$ & $1.62829\times10^{-19}$ & $1.38295\times10^{-22}$ & $1.17458\times10^{-25}$ \\
	\hline
	$0.6$ & $3.00672\times10^{-13}$  & $2.55368\times10^{16}$ & $2.16891\times10^{-19}$ & $1.84211\times10^{-22}$ & $1.56455\times10^{-25}$ \\
	\hline
	$0.7$ & $3.78205\times10^{-13}$  & $3.21219\times10^{-16}$ & $2.72819\times10^{-19}$ & $2.31713\times10^{-22}$ & $1.96799\times10^{-25}$ \\
	\hline
	$0.8$ & $4.55765\times10^{-13}$  & $3.87093\times10^{-16}$ & $3.28768\times10^{-19}$ & $2.79231\times10^{-22}$ & $2.37158\times10^{-25}$ \\
	\hline
	$0.9$ & $5.31029\times10^{-13}$  & $4.51017\times10^{-16}$ & $3.8306\times10^{-19}$ & $3.25343\times10^{-22}$ & $2.76322\times10^{-25}$ \\
	\hline
	$1$ & $6.01928\times10^{-13}$  & $5.11233\times10^{-16}$ & $4.34203\times10^{-19}$ & $3.6878\times10^{-22}$ & $3.13214\times10^{-25}$ \\
	\hline
\end{tabular}
\\
\\
It is observed that MDJM solution has less error than other iterative methods described above. Also we have used fewer terms of MDJM series than other methods to approximate the solution.
\\

\section{Conclusions}\label{concl}
The Boussinesq equation is an important class of PDEs arising in applied science. Various authors proposed different methods to solve these equations. In this article, we proposed a modification to DJM viz. MDJM and used it to solve few Boussinesq equations. It is observed that the MDJM is simpler iterative method than other iterative methods. Further, the solutions obtained by using this new method are good approximation to the exact solutions. This new method can be used to solve different nonlinear problems in a more efficient way.

\subsection{Acknowledgements}

S. Bhalekar acknowledges the Science and Engineering Research Board (SERB), New Delhi, India for the Research Grant (Ref. MTR/2017/000068) under Mathematical Research Impact Centric Support (MATRICS) Scheme


\section{References}
\begin{thebibliography}{99}
	
	\bibitem{BE} J. Boussinesq ,  Théorie de l’intumescence liquide appelée onde solitaire ou de translation se propageant dans un canal rectangulaire, Comptes Rendus Acad. Sci (Paris) 72 (1871), 755--759.
	
	\bibitem{lattices} Hirota, Ryogo, Exact N-soliton solutions of the wave equation of long waves in shallow-water and in nonlinear lattices, Journal of Mathematical Physics 14(7) (1973), 810--814.
	
	\bibitem{string} V. E. Zakharov, On stochastization of one-dimensional chains of nonlinear oscillators, Soviet Journal of Experimental and Theoretical Physics 38 (1974), 108.
	
	\bibitem{dielectric} S. K. Turitsyn, Nonstable solitons and sharp criteria for wave collapse, Physical Review E 47(1) (1993), R13.
	
	
	\bibitem{VIM}  Mo. Jiaqi, A variational iteration solving method for a class of generalized Boussinesq equations, Chin Phys Lett, 26(6) (2009), 060202.
	
	
	\bibitem{MVIM1} T. A. Abassy, M. A. El-Tawil, H. El-Zoheiry,  Modified variational iteration method for Boussinesq equation, Computers and Mathematics with Applications,  54(7) (2007), 955--965.
	

	
	\bibitem{MVIM2} E. Salehpour, H. Jafari, C. M. Khalique,  A modified variational iteration method for solving generalized Boussinesq equation and Lie’nard equation. International Journal of Physical Sciences, 6(23) (2011), 5406--5411.
	
	\bibitem{Mohyud} Mohyud-Din, Syed Tauseef, On Numerical Solutions of Two-Dimensional Boussinesq Equations by Using Adomian Decomposition and He's Homotopy Perturbation Method, J. Application and Applied Mathematics, 9466(2010) (1932) 1--11.
	
	\bibitem{ADM1} A. G. Bratsos, I. T. H. Famelis, D. P. Papadopoulos, On the solution of the Boussinesq equation using the Adomian decomposition method, (2004), 406-409.
	
	\bibitem{Malek} Abd-el-Malek, B. Mina , A. Nagwa Badran, Hossam S. Hassan, and Heba H. Abbas, New solutions for solving Boussinesq equation via potential symmetries method, Applied Mathematics and Computation, 251 (2015), 225-232.
	
	\bibitem{GEJJI1} V. Daftardar-Gejji, H. Jafari, An iterative method for solving non linear functional equations, J. Math. Anal. Appl. 316 (2006), 753--763.
	
	\bibitem{GEJJI3} V. Daftardar-Gejji, S. Bhalekar,  Solving fractional diffusion-wave equations using the New Iterative Method, Frac. Calc. Appl. Anal. 11 (2008), 193--202.
	
	\bibitem{GEJJI4} S. Bhalekar, V. Daftardar-Gejji,  New iterative method: Application to partial differential equations, Appl. Math. Comput. 203 (2008), 778--783.
	
	\bibitem{GEJJI5} V. Daftardar-Gejji, S. Bhalekar,  Solving fractional boundary value problems with Dirichlet boundary conditions, Comput. Math. Appl.  59 (2010), 1801--1809.
	
	\bibitem{GEJJI6} S. Bhalekar, V. Daftardar-Gejji,  Solving evolution equations using a new iterative method, Numer. Methods Partial Differential Equations 26 (2010), 906--916.
	
	\bibitem{GEJJI8} S. Bhalekar, V. Daftardar-Gejji,  Solving a system of nonlinear functional equations using revised new iterative method, Int. J. Computa. Math. Sci. 6 (2012), 127--131.
	
	\bibitem{Noor} M. A. Noor, K. I. Noor,  S. T. Mohyud-Din, A. Shabbir,  An iterative method with cubic convergence for nonlinear equations, Applied Mathematics and Computation, 183(2) (2006),  1249-1255.
	
	\bibitem{pb1} J. Patade, S. Bhalekar, Approximate analytical solutions of Newell-Whitehead-Segel equation using a new iterative method. World Journal of Modelling and Simulation, 11(2) (2015), 94-103.
	
	\bibitem{pb8} S. Bhalekar, J. Patade, An Analytical Solution of Fisher’s Equation Using Decomposition Method. American Journal of Computational and Applied Mathematics, 6(3) (2016), 123-127. 
	
	\bibitem{pb5} J. Patade, S. Bhalekar, On Analytical Solution of Ambartsumian Equation. Natl. Acad. Sci. Lett., 40(4) (2017),  291–293. 
	
		\bibitem{GEJJI7} S. Bhalekar, V. Daftardar-Gejji, Solving fractional order logistic equation using a new iterative method, Int. J. Differ. Equ. 2010 (2012), (Art. ID 975829).
	
	\bibitem{GEJJI9} S. Bhalekar, V. Daftardar-Gejji,  Numeric-analytic solutions of dynamical systems using a new iterative method, J. Appl. Nonlin. Dyn. 1 (2012), 141--158.
	
	\bibitem{pb4} S. Bhalekar, J. Patade, Analytical solutions of nonlinear equations with proportional delays. Appl. Comput. Math., 15(3) (2016),  331--345. 
		
	\bibitem{pb2}  S. Bhalekar, J. Patade, Series Solution of the Pantograph Equation and Its Properties. Fractal and Fractional, 1(1) (2017), 16. 	
	
	\bibitem{pb3}   S. Bhalekar, J. Patade,  Analytical Solution of Pantograph Equation with Incommensurate Delay. Physical Sciences Reviews, 2(9),  (2017)
	

	
	
	\bibitem{pb6} J. Patade, S. Bhalekar, A new numerical method based on Daftardar-Gejji and Jafari technique for solving differential Equations.  World J. Modell. Simul., 11  (2015)   256--271.
	
	\bibitem{pb7}	S. Bhalekar, J. Patade,   A Novel Numerical Method for Solving Volterra      
	Integro-Differential Equation.  Int. J. Appl. Comput. Math.,  6(1) (2020), 1-19.
	
	\bibitem{CONV} S. Bhalekar , V. Daftardar-Gejji,  Convergence of the new iterative method, Int. J. Differ. Equ. 2011 (2011).
	
	\bibitem{Wazwaz} A. M. Wazwaz, A reliable modification of Adomian decomposition method, Applied Mathematics and Computation 102(1) (1999), 77--86.
	
	\bibitem{Clarkson} P. A. Clarkson , M. D. Kruskal, New similarity reductions of the Boussinesq equation, Journal of Mathematical Physics 30(10) (1989), 2201--2213.
	
\end{thebibliography}
\end{document}